\documentclass[11pt]{amsart}
\usepackage{amsmath,amssymb,amsfonts,amsthm,amscd,indentfirst}

\usepackage[dvipdfm]{hyperref}

\numberwithin{equation}{section}

\newtheorem{prop}{Proposition}[section]
\newtheorem{theo}[prop]{Theorem}

\newtheorem{rema}[prop]{Remark}

\newcommand \xLongLeftRightArrow[2][]{ %
\ext@arrow 0055{\LongLeftRightArrowfill@ }{#1}{#2}}
\def\LongLeftRightArrowfill@ {%
\arrowfill@ \Leftarrow \Relbar \Rightarrow }

\newcommand \xlongleftrightarrow[2][]{
\ext@arrow 0055{\longleftrightarrowfill@}{#1}{#2}}
\def\longleftrightarrowfill@ {
\arrowfill@ \leftarrow \relbar \rightarrow }

\def\begeq{\begin{equation}}
\def\endeq{\end{equation}}

\def\Dint{\displaystyle\int}

\begin{document}
\title[Log Entropy functional]{ First variation of the Log Entropy functional along the Ricci flow}
\author{J\MakeLowercase{unfang} L\MakeLowercase{i}}
\address{Department of Mathematics\\
         McGill University\\
         Montreal, Quebec. H3A 2K6, Canada.}
\email{jli@math.mcgill.ca}
\thanks{Research of the author is
supported in part by a CRM fellowship.}
\date{} \maketitle

\begin{abstract}
In this note, we establish the first variation formula of the
adjusted log entropy functional $\mathcal Y_a$ introduced by Ye in
\cite{Y2}. As a direct consequence, we also obtain the monotonicity
of $\mathcal Y_a$ along the Ricci flow.
\end{abstract}

\setcounter {section} {1}
\bigskip

Various entropy functionals play crucial role in the singularity
analysis of Ricci flow. Let $(M^n, g(t))$ be a smooth family of
Riemannian metrics on a closed manifold $M^n$ and suppose $g(t)$ is
a solution of Hamilton's Ricci flow equation. In a recent
interesting paper \cite{Y2}, R. Ye introduced a new entropy
functional, the adjusted log entropy, as follows
\begin{equation}
\begin{array}{rll}
\mathcal Y_a(g,u,t)=-\Dint_Mu^2\ln u^2dvol + \frac{n}{2}\ln
\bigg(\Dint_M(|\nabla u|^2+\frac{R}{4}u^2)dvol+a \bigg) +4at,
\end{array}
\end{equation}
where the positive function $u\in W^{1,2}(M^n)$ satisfies
$\int_M(|\nabla u|^2+\frac{R}{4}u^2)dvol+a>0$, and $R$ denotes the
scalar curvature of the metric at time $t$.

The log entropy functional can be used to prove uniform logarithmic
Sobolev inequalities along the Ricci flow which also leads to
uniform Sobolev inequalities, see Ye's recent series of papers,
\cite{Y1}, \cite{Y2}, etc, and Zhang \cite{Z}. This new entropy
functional of Ye shares a similar important feature with Perelman's
entropy functionals. Namely, it is nondecreasing under the following
coupled system of Ricci flow,
\begin{equation}\label{coupled system}
\left\{
\begin{array}{rll}
\frac{\partial}{\partial t}g_{ij}=&-2R_{ij}\\
\frac{\partial}{\partial t} u=& -\Delta u-\frac{|\nabla
u|^2}{u}+\frac{R}{2}u.
\end{array}
\right.
\end{equation}

The first evolution equation is the Ricci flow equation. The second
equation ensures $\Dint_M u^2d\mu_g=1$ to be preserved by the Ricci
flow. Notice that, if we define $u=e^{-\frac{f}{2}}$, then the
second equation is equivalent to Perelman's equation
$\frac{\partial}{\partial t} f= -\Delta f+{|\nabla f|^2}-R$ which
instead preserves $\Dint_M e^{-f}d\mu_g=1$.

 The following statement is obtained by Ye,  see
{\rm {\bf Theorem 3.1}} in \cite{Y2}.\\

 {\em Assume
$a>-\lambda_0(g)$. Then $\mathcal Y_a\equiv\mathcal
Y_a(g(t),u(t),t)$ is nondecreasing. Indeed, we have
\begin{equation}\label{Ye equ1}
\begin{array}{rll}
\frac{d}{dt} \mathcal
Y_a\ge\displaystyle\frac{n}{4\omega}\Dint_M|Ric-2\frac{\nabla^2u}{u}+2\frac{\nabla
u\otimes\nabla u}{u^2} -\frac{4\omega}{n}{g}|^2u^2d\mu_g,
\end{array}
\end{equation}
where $\omega = \omega(t) = a + \Dint_M(|\nabla
u|^2+\frac{R}{4}u^2)dvol\bigg|_t$,
which is positive.\\

}

 Ye used a minimizing procedure and the monotonicity formula
of $\mathcal W$-entropy of Perelman to show the monotonicity formula
of (\ref{Ye equ1}),
see {\bf Lemma 4.1}, {\bf4.2} in \cite{Y2}.\\

One question may be interesting is : what is the precise first
variation of the functional $\mathcal Y_a$? In this short note,
instead of giving an inequality for the first variation, we derive a
formula of the first variation itself for $\mathcal Y_a$.
Consequently, we also obtain the monotonicity of $\mathcal Y_a$
along the Ricci flow. The proof is a direct approach without
appealing to the minimizing procedure and the monotonicity of
$\mathcal W$ functional. In Remark \ref{Remark on Ye}, we will show
that the first variation formula (\ref{equ1}) we obtained is
equivalent to the righthand side of Ye's inequality (\ref{Ye equ1}),
see Remark \ref{Remark on Ye}.

We adapt the notations in \cite{Y2}. If we let $u=e^{-\frac{f}{2}}$,
then $\mathcal F=\int_M (R+|\nabla f|^2)e^{-f}d\mu_g=4\int_M
(|\nabla u|^2+\frac{1}{4}Ru^2)d\mu_g$ which implies
$4\omega=4a+\mathcal F$. We note that $\mathcal F$ is one of the
entropy functionals introduced by Perelman \cite{P}. Now we
introduce the main theorem of this paper.


\begin{theo}
Assume $a>-\lambda_0(g)\ (=-\frac{1}{4}\lambda_0(\mathcal F))$. Then
$\mathcal Y_a\equiv\mathcal Y_a(g(t),u(t),t)$ is nondecreasing.
Indeed, we have
\begin{equation}\label{equ1}
\begin{array}{rll}
\frac{d}{dt} \mathcal
Y_a=&\displaystyle\frac{n}{4\omega}\int_M|Ric-2\frac{\nabla^2u}{u}+2\frac{\nabla
u\otimes\nabla u}{u^2}-\frac{4(\omega-a)}{n}g|^2u^2d\mu_g
+\frac{4a^2}{\omega},
\end{array}
\end{equation}
where $\omega = a + \Dint_M(|\nabla
u|^2+\frac{R}{4}u^2)d\mu_g\bigg|_t$, which is positive. The
monotonicity is strict, unless the manifold is a gradient shrinking
soliton with positive first eigenvalue $\lambda_0(g)>0$ and also
$a=0$.
\end{theo}

\begin{proof}
By the notations of $u$ and $f$, we observe that $\mathcal
Y_a=-\mathcal S + \frac{n}{2}\ln\frac{1}{4}(\mathcal F+4a)+4at $,
where $\mathcal S=-\int_Mfe^{-f}d\mu_g$ is the differential Shannon
Entropy. Along the coupled system of Ricci flow, we have
\begin{equation}\label{}
\begin{array}{rll}
\frac{d}{dt} \mathcal Y_a=&4a-\mathcal
F+\frac{n}{2(\mathcal F+4a)}\frac{d}{dt}\mathcal F\\
=& \frac{n}{2(\mathcal F+4a)}\big[\frac{d}{dt}\mathcal F+\frac{2}{n}(4a-\mathcal F)(4a+\mathcal F)\big]\\
=& \frac{n}{(\mathcal F+4a)}\big[\Dint_M|R_{ij}+\nabla_i\nabla_j
f|^2e^{-f}d\mu_g
+\frac{1}{n}(4a-\mathcal F)(4a+\mathcal F)\big]\\
=&\frac{n}{4\omega}\big[\Dint_M|R_{ij}+\nabla_i\nabla_j
f-\frac{\mathcal F}{n}g_{ij}|^2e^{-f}d\mu_g
+\frac{(4a)^2}{n}\big]\ge 0.\\
\end{array}
\end{equation}
Change of variables from $f$ to $u$ completes the proof.
\end{proof}
In the above, we have used the properties that $ \frac{d}{dt}
\mathcal S=\mathcal F$ and $\frac{d}{dt} \mathcal F=
2\int|R_{ij}+\nabla_i\nabla_j f|^2e^{-f}d\mu_g$ under the coupled
system (\ref{coupled system}) which can be proved by direct
computations. Related references can be found in \cite{CLN},
\cite{BM}, or original paper of Perelman \cite{P} and \cite{N2},
\cite{L1}. Recall that, $\lambda_0(g)$ denotes the first eigenvalue
of $-\Delta +\frac{R}{4}$ in \cite{Y2}, i.e.
\[
\lambda_0(g)=\inf\int_M(|\nabla u|^2+\frac{R}{4}u^2)d\mu_g,
\]
where the infimum is taken over all $u$ satisfying $\int_M u^2d\mu_g
= 1$. Hence, the second property ensures that $\lambda_0(g(t))$ is
nondecreasing and $\int_M(|\nabla u|^2+\frac{R}{4}u^2)d\mu_g+a>0 $
for all time $t$.

 Instead of working on $\mathcal W$-functional, dealing
directly with $\mathcal S$ and $\mathcal F$-functionals makes the
proof rather elementary. Below a few remarks are in order.

\begin{rema}
From the first variation formula, we know precisely when the
monotonicity is strict. Only  when $a=0$ the monotonicity can be
non-strict. In this case, we know the manifold must be a shrinking
gradient Ricci soliton with $\lambda_0(g)>0$.
\end{rema}

\begin{rema}\label{Remark on Ye}
Simple computations yields that (\ref{equ1}) we obtained is
equivalent to the righthand side of Ye's inequality (\ref{Ye equ1}).
\begin{equation}\label{}
\begin{array}{rll}
\frac{d}{dt} \mathcal Y_a
=&\frac{n}{4\omega}\bigg[\Dint_M|R_{ij}+\nabla_i\nabla_j
f-\frac{4(\omega-a)}{n}g_{ij}|^2e^{-f}d\mu_g +\frac{(4a)^2}{n}\bigg]\\
=&\frac{n}{4\omega}\bigg[\Dint_M|R_{ij}+\nabla_i\nabla_j
f-\frac{4\omega}{n}g_{ij}|^2e^{-f}d\mu_g \\
&+2\Dint_Mg^{ij}\big(R_{ij}+\nabla_i\nabla_j
f-\frac{4\omega}{n}g_{ij}\big)\frac{4a}{n}e^{-f}d\mu_g + \frac{(4a)^2}{n}+\frac{(4a)^2}{n}\bigg]\\
=&\frac{n}{4\omega}\bigg[\Dint_M|R_{ij}+\nabla_i\nabla_j
f-\frac{4\omega}{n}g_{ij}|^2e^{-f}d\mu_g \bigg].\\
\end{array}
\end{equation}

The splitting sum we had in (\ref{equ1}) clearly shows that when
$a\neq 0$, the monotonicity is strict.
\end{rema}

\begin{rema}
There are intensive study on various entropy functionals which are
related to or motivated by Perelman's entropy functionals in
\cite{P}. For example, entropy functionals for linear heat equations
in \cite{N1}, expanding $\mathcal W$-entropy in \cite{FIN}, Log
entropy in \cite{Y2}, entropy on fiber bundles in \cite{Lott},
entropy on an extended Ricci flow system in \cite{List}, and various
generalized entropy functionals in \cite{L2} and \cite{KZ}. Also we
note that there is an interesting entropy functional on the
evolution equation for $p$-harmonic functions appeared in \cite{KN}
recently. All these entropies share one common feature : they have
monotonicity properties under geometric evolution equations and the
monotonicity is strict unless the metric is a soliton.
\end{rema}

\bigskip
{\em Acknowledgement.} I would like to thank professor B. Chow for
discussions.

\end{document}